\documentclass[11pt]{article}
\usepackage[a4paper, margin=0.8in, bmargin=0.8in]{geometry}
\usepackage[mathscr]{euscript}
\usepackage{amsmath,amsthm}
\usepackage{amsfonts,amssymb}
\usepackage{graphicx}
\usepackage{breqn}
\usepackage{mathrsfs}
\usepackage{amsmath}
\usepackage{float}
\usepackage{caption}
\usepackage{chngcntr}
\usepackage{tikz-qtree}
\usetikzlibrary{arrows.meta}
\usepackage{subcaption}
\usepackage{epsfig,amsmath,amsthm,latexsym,amssymb,amsgen,graphicx}
\usepackage{verbatim}
\usepackage{cite}
\newtheorem{thm}{Theorem}[section]

\newtheorem{cor}[thm]{Corollary}

\newtheorem{prop}[thm]{Proposition}

\newtheorem{defn}[thm]{Definition}
\newtheorem{rem}[thm]{Remark}

\newtheorem{exam}[thm]{Example}
\begin{document}
\title {\bf{Strong Domination Index in Fuzzy Graphs}}
\author { \small Kavya. R. Nair \footnote{Corresponding author. Kavya. R. Nair, Department of Mathematics, National Institute of Technology, Calicut, Kerala, India. Email: kavyarnair@gmail.com} and M. S. Sunitha \footnote{sunitha@nitc.ac.in} \\ \small Department of Mathematics, National Institute of Technology, Calicut, Kerala, India-673601}
\date{}
\maketitle
\hrule
\begin{abstract}
Topological indices play a vital role in the area of graph theory and fuzzy graph (FG) theory. It has wide applications in the areas such as chemical graph theory, mathematical chemistry, etc. Topological indices produce a numerical parameter associated with a graph. Numerous topological indices are studied due to its applications in various fields. In this article a novel idea of domination index in a FG is defined using weight of strong edges. The strong domination degree (SDD) of a vertex $u$ is defined using the weight of minimal strong dominating set (MSDS) containing $u$. Idea of upper strong domination number, strong irredundance number, strong upper irredundance number, strong independent domination number, and strong independence number are explained and illustrated subsequently. Strong domination index (SDI) of a FG is defined using the SDD of each vertex. The concept is applied on various FGs like complete FG, complete bipartite and r-partite FG, fuzzy tree, fuzzy cycle and fuzzy stars. Inequalities involving the SDD and SDI are obtained. The union and join of FG is also considered in the study. Applications for SDD of a vertex is provided in later sections. An algorithm to obtain a MSDS containing a particular vertex is also discussed in the article.
\end{abstract}
\hrule
\section{Introduction}
One of the most efficient and precise methods to express information and the relationships between different entities is through a graph. A weighted graph shows the strength of the relations between the vertices. Due to its use in numerous fields, domination in graphs is one of the most important topics of research. The research on domination was pioneered by Ore \cite{ore1962} and Berge \cite{berge1962} in 1962, and it was further developed by Hedetniemi \cite{haynes} and Cockayne \cite{bollobas, cockayne1, cockayne2, cockayne3}. Fuzzy graphs (FGs) were developed as a result of the element's uncertain or ambiguous character and the way they are related to one another. In 1975, Rosenfeld \cite{rosenfeld} established the study of fuzzy graphs and the fuzzy equivalents of graph notions such as trees, paths, cycles, connectedness, etc. FG models are better than graph models since they consider the ambiguity and uncertainty in the entities under consideration. Mordeson et al. \cite{book2} discussed the connectivity concepts and fuzzy end nodes in FGs. The idea of $\alpha-, \beta-, \delta- $ edges was introduced by Mathew and Sunitha \cite{sm}. Domination in FGs using the weight of strong arcs was developed in \cite{manjusha}. Mordeson and Peng \cite{mordeson} discussed operations in FGs. A topological index, also known as a molecular structure descriptor, is a number that relates to the composition of chemicals and is used to compare the chemical structure to different physical attributes, chemical reactivity, or biological activity. In addition to being beneficial in several chemical fields, molecular-graph-based topological indices have also been shown to be helpful in computational biology and computational linguistics. Presently, there are countless similar indices in the literature. 
Distance-based topological indices and degree-based topological indices are two of the main types of topological indices. Degree-based topological indices are extremely important. The first distance-based index was the Wiener index W(G), which was developed in 1974 by the chemist Wiener \cite{W}. Later, a number of indices developed as a consequence of their application in chemistry; some of them are the Zagreb indices, Randic, Harmonic, and Gutman and Schultz indices. Wiener and connectivity indices were introduced into fuzzy graphs by Binu et al. \cite{binu}. Other topological indices, including the modified and hyper-Wiener, Gutman, Schultz, Zagreb, Harmonic, and Randic indices in FGs, are investigated in \cite{T6}. Ample work on indices is done in \cite{T1, T2, T3, T4}. Most of the works on topological indices are extended to other fuzzy graphs like bipolar FGS, fuzzy incidence graphs and intuitionistic FGs. This work is inspired from the domination index on graphs defined by Ahmed et al. \cite{T5}. In the article the domination degree of a vertex $v$ is defined as the number of minimal dominating sets containing $v.$ Motivated from this definition, this article discusses the domination degree and domination index in FGs using weight of strong arcs. 
The concept of strong domination in FGs uses the weight of strong edges and hence the least value for domination number is obtained in this concept. Also, to define the SDD of a vertex the concept of MSDS is used rather than minimum SDS, since it is not always possible to find a minimum SDS containing a particular vertex. But it is always possible to find a MSDS containing a particular vertex, which is proved in Theorem \ref{37}. \\
The article is structured as follows: Section 2 deals with the preliminaries. Section 3 begins with the theorem that ascertains that it is always possible to find a MSDS containing a particular vertex in a connected SFG. The remainder of the section defines the SDD of a vertex, a strong domination regular FG, upper strong domination number, strong irredundance number, strong upper irredundance number, strong independent domination number, and strong independence number and establishes some inequalities and results related to the definitions. Section 4 contains the definition of strong domination index (SDI). The SDI of certain FGs are studied and illustrated in the section. The SDI is studied in fuzzy trees, fuzzy cycles, union and join of FGs. Applications for the SDD of a vertex is provided in Section 5. Section 6 deals with an algorithm to find a MSDS containing a particular vertex.
\section{Preliminaries}
Throughout the article , the symbols $\wedge$ and $\vee$ represent minimum and maximum operators respectively. The basic concepts of FGs are referred from \cite{zadeh, book, book2, book4, sm, bhutani, anjali, binu, nagoor2008, nagoor2006, manjusha, mordeson, book3  }.
A fuzzy graph (FG), $\mathcal{X}=(\mathcal{V}, \varrho, \upsilon )$ is a triple, such that $\mathcal{V}\neq \phi,$ $\varrho: \mathcal{V}\rightarrow [0,1]$, $\upsilon: \mathcal{V}\times \mathcal{V}\rightarrow [0,1]$ and $\upsilon(ab)\leq \varrho(a)\wedge \varrho(b)$ $\forall a,b\in \mathcal{V}.$ Here $\varrho $ and $\upsilon$ are the fuzzy vertex and edge sets of $\mathcal{X} $ respectively and $\varrho^*=\{a\in \mathcal{V}: \varrho(a)>0\}$, $\upsilon^*=\{ab\in \mathcal{V}\times \mathcal{V}: \upsilon(ab)>0\}$. Hence, $\mathcal{X}^*=(\rho^*,\upsilon^*)$ is the underlying graph of $\mathcal{X}.$ Order $\mathfrak{p}$ and size $\mathfrak{q}$ are respectively defined as $\mathfrak{p}=\sum\limits_{a\in \varrho^*} \varrho(a),$ $\mathfrak{q}= \sum\limits_{ab\in \upsilon^*}\upsilon(ab).$ The FG, $\mathcal{Y}=(\mathcal{V}, \varepsilon, \tau)$ is a partial fuzzy subgraph of $\mathcal{X} $ if $\varepsilon \subseteq \varrho $ and $\tau\subseteq \upsilon.$ Also, $\mathcal{Y}$ is spanning fuzzy subgraph if $\varepsilon=\varrho.$ A sequence of distinct vertices $a_1,a_2,...,a_n$ such that $\upsilon(a_i,a_{i+1})>0$, $i=1,2,...n-1$ is called a path $P$ from $a_1 $ to $a_n$. The minimum of weight of edges in a path $P$ is called the strength of the path, denoted as $S(P)$. The maximum of strengths of all paths from $a$ to $b$ is called strength of connectedness between $a$ and $b$, denoted as $CONN_{\mathcal{X}}(a,b)$ or $\upsilon^{\infty}(a,b).$ Two vertices are connected if they are joined by a path. If every pair of vertices are connected then the FG, $\mathcal{X}$ is said to be connected.
Consider two FGs, $\mathcal{X}=(\mathcal{V},\varrho,\upsilon)$ and $\mathcal{X}'=(\mathcal{V}',\varrho',\upsilon').$ A bijective map $\Theta: \mathcal{V}\rightarrow \mathcal{V}'$ such that $\varrho(a)=\varrho'(\Theta(a))$ and $\upsilon(a,b)=\upsilon'(\Theta(a),\Theta(b))$ $\forall a,b\in \mathcal{V}$, is called fuzzy isomorphism from $\mathcal{X}$ to $\mathcal{X}'$. 
If $\upsilon(ab)>\upsilon^{'\infty}(a,b)$, where $\upsilon^{'\infty}(a,b)$ is the strength of connectedness between $a$ and $b$ when $ab$ is deleted from $\mathcal{X}, $ then $ab$ is called an $\alpha-$ strong edge. Also, $ab $ is a $\beta-$ strong edge if $\upsilon(ab)=\upsilon^{'\infty}(a,b)$. A $\delta-$ edge $ab$, is such that $\upsilon(ab)<\upsilon^{'\infty}(a,b)$. Here $\alpha $ and $\beta-$ strong edges are called strong edges. A FG, $\mathcal{X}$ is called strong fuzzy graph (SFG) if all edges in $\mathcal{X}$ are strong. A path having only strong edges is called a strong path.
A strong path $P$ from $a $ to $b$ is an $a-b$ geodesic if there is no shorter strong path from $a$ to $b$. The sum of weights of the edges in a geodesic is called the weight of the geodesic. Now, $\mathcal{X}$ is a cycle if $(\varrho^*, \upsilon^*)$ is a cycle. If in addition, if $\mathcal{X}$ has the property that there is at least two edges with the least weight then $\mathcal{X}$ is a fuzzy cycle (FC). A FC is $\beta-$ saturated if each vertex has at least one $\beta-$ edge incident to it. The FG, $\mathcal{X}=(\varrho,\upsilon)$ is a tree if $(\varrho^*,\upsilon^*)$ is a tree. If $\mathcal{X}$ has a fuzzy spanning subgraph $\mathcal{F}=(\varrho,\tau)$ which is tree and $\upsilon(ab)<\tau^{\infty}(a,b)$, $\forall ab\in \upsilon^*\setminus \tau^*$, then $\mathcal{X}$ is a fuzzy tree (FT). Now, let $\mathcal{X}$ be a connected FG. A maximum spanning tree (MST) is a spanning subgraph $\mathscr{T}=(\varrho, \tau)$, such that $\mathscr{T}$ is a tree and $S(P)$ is $\upsilon^{\infty}(a,b)$, where $P$ is the unique path between $a$ and $b$ $\forall a,b\in \varrho^*.$ A FG is complete FG (CFG), if $\upsilon(ab)=\varrho(a)\wedge\varrho(b), \forall a,b\in\varrho^*.$ For a bipartite FG, the vertex set can be partitioned into two sets $\mathcal{V}_1$ and $\mathcal{V}_2$ that are non- empty and $\upsilon(ab)=0$ if $a,b\in \mathcal{V}_1$ or $a,b\in \mathcal{V}_2$. Also, if $\upsilon(ab)=\varrho(a)\wedge\varrho(b), \forall a\in \mathcal{V}_1$ and $b\in \mathcal{V}_2$, then $\mathcal{X}$ is a complete bipartite FG (CBFG), denoted by $\mathcal{K}_{\varrho_1,\varrho_2}$, where $\rho_1$ and $\rho_2$ are restrictions of $\rho$ to $\mathcal{V}_1$ and $\mathcal{V}_2$ respectively. If the unique MST of a fuzzy tree, $\mathcal{X}$ is a star, then $\mathcal{X}$ is called a fuzzy star (FS), denoted as $FS_n, $ where $|\varrho^*|=n+1.$ \\
Let $\mathcal{X}$ be a FG, two vertices $a$ and $b$ are strong neighbors if $ab$ is a strong edge. Here the closed strong neighborhood is $N_s[a]=\{a\}\cup\{b\in\varrho^*: ab \text{ is a strong edge}\}$. The private neighborhood of vertex $a$ with respect to set $\mathcal{S}\subseteq \mathcal{V}$, denoted as $P_{N_s}[a,\mathcal{S}]=N_s[a]\setminus \bigcup \limits _{b\in\varrho^*\setminus a}N_s[b]$. If $P_{N_s}[a,\mathcal{S}]\neq \phi, \forall a\in \mathcal{S}$, then $\mathcal{S} $ is called a fuzzy irredundant set (FIRS). A FIRS $\mathcal{S}$, is a maximal FIRS (MFIRS), if for every $a\in \mathcal{V}\setminus \mathcal{S} $, $\mathcal{S}\cup \{a\}$ is not a FIRS. If the edge between two vertices $a$ and $b$ is not strong, then $a$ and $b$ are said to be fuzzy independent. A set $\mathcal{M}\subseteq \mathcal{V}$ is fuzzy independent set (FIS), if every pair of vertices in $\mathcal{M}$ is fuzzy independent. A vertex $a$ strongly dominates itself and all its strong neighbors. A set $\mathcal{D}\subseteq \mathcal{V} $ is called a strong dominating set (SDS) if $\forall b \in \mathcal{V} \setminus \mathcal{D}$,  $\exists a\in \mathcal{D}$ such that $a$ dominates $b.$ A SDS, $\mathcal{D}$ is minimal (MSDS) if $\mathcal{D}\setminus a, a\in \mathcal{D}$ is not a SDS. The weight of a set $\mathcal{S}\subseteq \mathcal{V} $ is $W(\mathcal{S})= \sum \limits_{a\in \mathcal{S}} \upsilon(ab)$, where $b\in N_s[a]$ and $ab $ is a strong edge with least weight incident at $a$. The minimum weight of SDS in $\mathcal{X} $ is called the strong domination number (SDN), $\gamma_s(\mathcal{X})$ or $\gamma_s$. A SDS with weight $\gamma_s$ is called a minimum SDS. An independent fuzzy dominating set (IFDS) is a set which is both a SDS and a FIS. \\
The Wiener index (WI) of $\mathcal{X}$ is defined as:
$$WI(\mathcal{X})=\sum \limits _{a,b\in \varrho^*} \varrho(a) \varrho(b) d_s(a,b),$$ where $d_s(a,b)$ is the minimum of the weights of geodesics from $a$ to $b.$ 
\\
Let $\mathcal{X}_1=(\mathcal{V},\varrho_1,\upsilon_1)$ and $\mathcal{X}_2=(\mathcal{V},\varrho_2,\upsilon_2)$ be two FGs with $\mathcal{X}_1^*=(\mathcal{V}_1,E_1)$ and $\mathcal{X}_2^*=(\mathcal{V}_2,E_2)$. Let $\mathcal{X}=(\mathcal{V},\varrho,\upsilon)$ be the union of $\mathcal{X}_1$ and $\mathcal{X}_2$ with edge set $E=E_1\cup E_2$. Here $\mathcal{V}=\mathcal{V}_1\cup \mathcal{V}_2$ and $\mathcal{V}_1\cap \mathcal{V}_2= \phi. $ Also $\varrho=\varrho_1\cup\varrho_2$ and $\upsilon=\upsilon_1\cup \upsilon_2$ are defined as follows:
$$\varrho(a)=\begin{cases}
\varrho_1(a)\, \text{ if }\, a\in \mathcal{V}_1\\
\varrho_2(a)\, \text{ if }\, a\in \mathcal{V}_2
\end{cases}
$$
$$\upsilon(ab)=\begin{cases}
\upsilon_1(ab)\, \text{ if }\, ab\in E_1\\
\upsilon_2(ab)\, \text{ if }\, ab\in E_2
\end{cases}
$$
Let $\mathcal{X}=(\mathcal{V},\varrho,\upsilon)$ be the join of $\mathcal{X}_1$ and $\mathcal{X}_2$ with edge set $E_1\cup E_2\cup E$, where $E=\{ab: a\in \mathcal{V}_1 \text{ and } b\in \mathcal{V}_2\}$. Here $\mathcal{V}=\mathcal{V}_1\cup \mathcal{V}_2$ and $\mathcal{V}_1\cap \mathcal{V}_2=\phi$. Also $\varrho=\varrho_1+\varrho_2$ and $\upsilon=\upsilon_1+\upsilon_2$ are defined as follows:
$$\varrho(a)=\begin{cases}
    \varrho_1(a) \text{ if } a\in \mathcal{V}_1\\
    \varrho_2(b) \text{ in } a\in \mathcal{V}_2
\end{cases}
$$
$$\upsilon(ab)=\begin{cases}
\upsilon_1(ab)\, \text{ if }\, ab\in E_1\\
\upsilon_2(ab)\, \text{ if }\, ab\in E_2\\
\varrho(a)\wedge\varrho(b) \text{ if } a\in \mathcal{V}_1 \text{ and } b\in \mathcal{V}_2
\end{cases}
$$
\begin{thm}\cite{nagoor2006}\label{56}
 Every minimum SDS of a FG is a minimal SDS.  
\end{thm}
\begin{thm}\cite{nagoor2006}\label{57}
Let $\mathcal{X}$ be a FG without isolated vertices. If $\mathcal{D}$ is a minimal SDS then $\mathcal{V}\setminus \mathcal{D}$ is a SDS.    
\end{thm}
\section{Domination index in fuzzy graphs}
A domination number is defined using a minimum dominating set. In graph, domination number is the number of vertices in a minimum dominating set. And in FGs it is the weight of a minimum dominating set. The weight of a dominating set can be defined in terms of weight of vertices or edges. But it is not always possible to find a minimum dominating set containing a particular vertex. Suppose $u$ is a vertex. Then it is not always possible to find a minimum dominating set containing $u$ in a FG $\mathcal{X}$. On the contrary it is always possible to find a minimal dominating set containing vertex $u$, $\forall u\in \mathcal{V}$. This result is proved in Theorem \ref{37}. Hence the strong domination degree (SDD) [Definition \ref{1}] of a vertex $u$ is defined using the weight of a MSDS containing $u$. The section discusses and illustrates the concept of SDD of a vertex. The notion of strong domination regular FG, upper strong domination number, strong irredundance number, strong upper irredundance number, strong independent domination number, and strong independence number are defined subsequently and inequalities involving these parameters are also obtained in this section.
\begin{thm}\label{37}
 Let $\mathcal{X}$ be a connected SFG and $u$ be a vertex in $\mathcal{X}$. Then it is always possible to find a MSDS containing $u$ in $\mathcal{X}$. 
\end{thm}
\begin{proof}
    Consider FG, $\mathcal{X}$ and a vertex $u\in \mathcal{V}(\mathcal{X})$. If $u\in D'$ such that $D'$ is a minimum SDS, then by Theorem \ref{56}, $D'$ is the required MSDS. Suppose $u$ does not belong to any minimum dominating set. Then for any minimum SDS $D$, $u\notin D$. Since $D$ is minimum SDS $D$ is minimal. Hence by Theorem \ref{57}, $\mathcal{V}\setminus D$ is a SDS containing $u.$ It is not necessary that $\mathcal{V}\setminus D$ is minimal. But it is always possible to construct a minimal dominating set containing $u$ from a dominating set containing $u$. In the following proof a MSDS containing $u$ is constructed from a SDS containing $u$.\\
    Let $\mathcal{S} $ be a SDS containing $u$. Suppose $\mathcal{S} $ is not minimal. Then $\exists$ a vertex $v$ such that $\mathcal{S} \setminus \{v\}$ is a dominating set. If $v\neq u$ and $\mathcal{S} \setminus \{v\}$ is MSDS then $\mathcal{S} \setminus \{v\}$ is the required MSDS. If not repeat the process and in each stage if $u$ is not the deleted vertex then finally a set containing $u$ is obtained which is a MSDS.\\
    Now suppose that $\mathcal{S} \setminus \{u\}$ is a SDS, i.e, $P_N[u, \mathcal{S}]=\phi.$ That means $u$ and all its neighbors are dominated by vertices in $\mathcal{S}$. \\
    \textbf{Case 1:} $N_s[u]$ is dominated by the same vertex say $v$ in $\mathcal{S}$.
    Suppose that $N_s[u]=\{u,u_1,u_2,...,u_n\}$. Also let $N_s[v]\setminus \{u,u_1,...u_n\}=\{v_1,v_2,...,v_m,v_1',v_2',...,v_k'\}$ where $v_1,v_2,..,v_m$ are end vertices. And there will be connected components containing $v$ as a vertex. The other $k$ vertices are neighbors of $v$ which are not end vertices that belong to a connected component containing $v$. Also among the $k$ neighbors of $v$ some vertices may be adjacent to each other. Now consider $\mathcal{S}$. Delete $v$ from $\mathcal{S}$ and add the following vertices to $\mathcal{S}$:
    \begin{enumerate}
        \item All the end vertices that are neighbors of $v$, i.e, $\{v_1,v_2,...,v_n\}.$
        \item If two vertices in $\{v_1',v_2',...,v_k'\}$ are neighbors then add any one of them. 
        \item If while considering two vertices $v_i'$ and $v_j'$ in $\{v_1',v_2',...,v_k'\}$ are not neighbors, then add both of them.
    \end{enumerate}
    Let $\mathcal{S}'$ be the newly obtained set. Now $P_N[u, \mathcal{S}']$ contains at least $u$ since the vertex dominating $u$ is deleted from $\mathcal{S}$. Hence, $P_N[u, \mathcal{S}']\neq \phi.$ If on addition of any of the vertices mentioned in $1,2$ and $3$ disturbs the minimality condition then delete the vertex with empty private neighborhood from $\mathcal{S}'$. And the resulting set will be a MSDS containing $u.$\\
    \textbf{Case 2:} $N_s[u]$ is dominated by distinct vertices in $\mathcal{S}$.\\
    Suppose $u$ is dominated by vertex $v$ in $\mathcal{S}$ and the neighbors of $u$ are dominated by $x_1,x_2,..x_{k'}$ in $\mathcal{S}$. Then repeat the same steps $1,2$ and $3$ with the neighbors of $v,x_1,x_2,...x_{k'}.$ and delete the vertex that disturbs the minimality condition. 
    The case when $u$ and its neighbors being dominated by multiple vertices of $\mathcal{S}$ can be proved similarly. \\
    Hence a MSDS containing a particular vertex $u$ is obtained from a SDS containing $u.$ Therefore, it is always possible to find a MSDS containing a particular vertex in a FG $\mathcal{X}$.   
\end{proof}
In Theorem \ref{37} a strong and connected FG $\mathcal{X}$ is considered. The case of non strong FG can be proved in a similar manner since the $\delta $ edges are neglected while considering strong domination. Also proof can be applied to each connected component in a disconnected FG.
The illustration of Theorem \ref{37} is provided in Example \ref{38}.

\vspace{0.8cm}
\begin{exam}\label{38}
Consider the FG in Figure \ref{39}. Here the marked vertices form a SDS containing $u$ which is not minimal. Here the private neighborhood of $u$ is empty. Now, the neighbors $\{u_1,u_2,u_3,u_4\}$ of $u$ and $u$ are dominated by $v,x_1,x_2$. Perform steps $1,2\text{ and }3$ and the set obtained is marked in Figure \ref{40}. Here the newly obtained SDS is not minimal. The vertices $y_1,y_2,y_3,\text{ and }y_4$ have empty private neighborhood in FG shown in Figure \ref{40} Hence delete $y_1,y_2,y_3,\text{ and }y_4$. Finally a SDS marked in Figure \ref{41} is obtained which is minimal and containing $u.$
\begin{figure}[h]
    \centering
    \includegraphics[width=12cm,height=6cm]{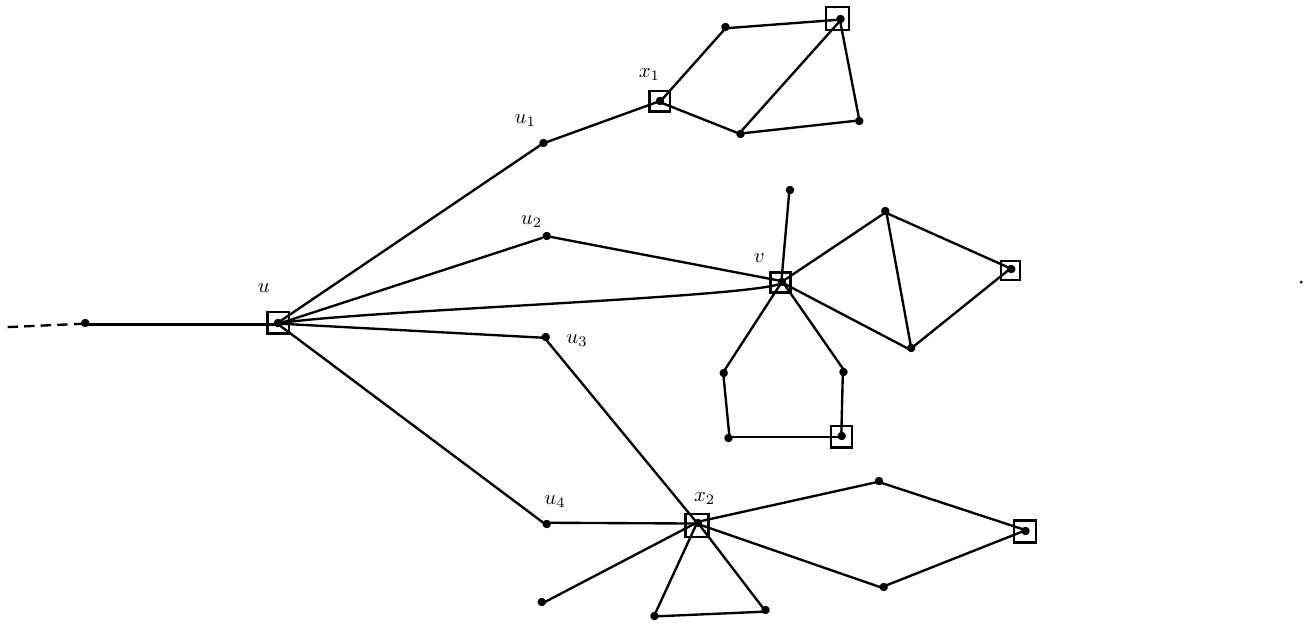}
    \caption{FG}
    \label{39}
\end{figure}
\begin{figure}[h]
    \centering
    \includegraphics[width=12cm,height=6cm]{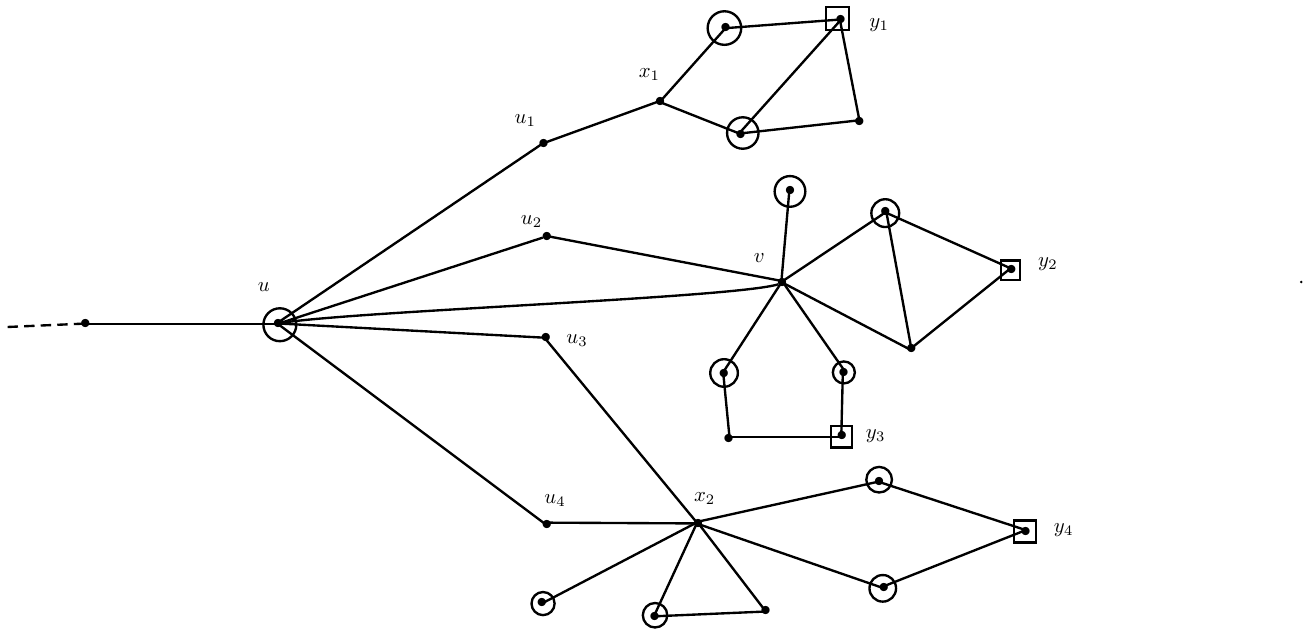}
    \caption{FG}
    \label{40}
\end{figure} 
\begin{figure}[h!]
    \centering
    \includegraphics[width=12cm,height=6cm]{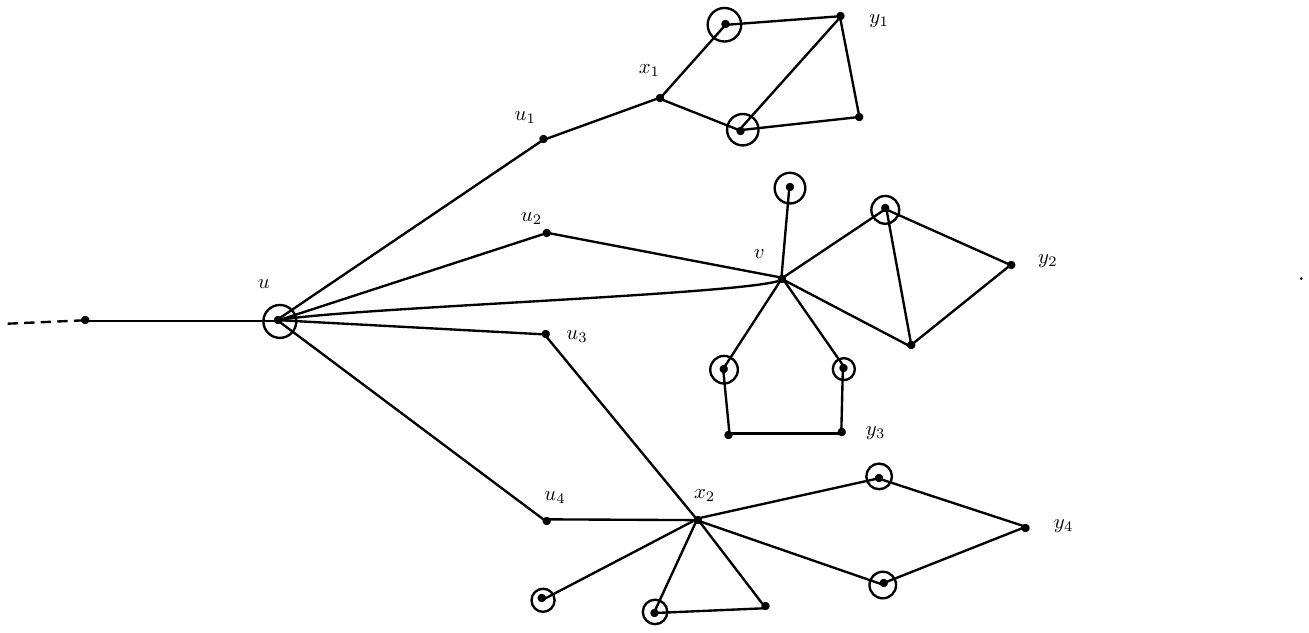}
    \caption{FG}
    \label{41}
\end{figure}
\end{exam}
\newpage
\begin{defn}\label{1}
Let $u$ be a vertex of FG, $\mathcal{X}$. Then the strong domination degree (SDD) of $u$ in $\mathcal{X}$ is the minimum of weight of MSDSs containing $u$, and is denoted as $_{sd}d_{\mathcal{X}}(u)$ or simply $_{sd}d(u)$, i.e,
$$_{sd}d(u)= min\, \{W(D): \, D \text{ is the MSDS containing } u\}$$
\end{defn}
Accordingly, the minimum and maximum SDD of a fuzzy graph are given by,
$$_{sp}\delta(G)= min\,\{_{sd}d(u):\, u\in V(\mathcal{X})\}$$
$$_{sp}\Delta(G)= max\,\{_{sd}d(u):\, u\in V(\mathcal{X})\}$$
The SDD of a vertex, minimum SDD of a FG and maximum SDD of a FG are explained using Example \ref{3}.
\begin{figure}[h]
    \centering
    \includegraphics[width=4cm,height=4cm]{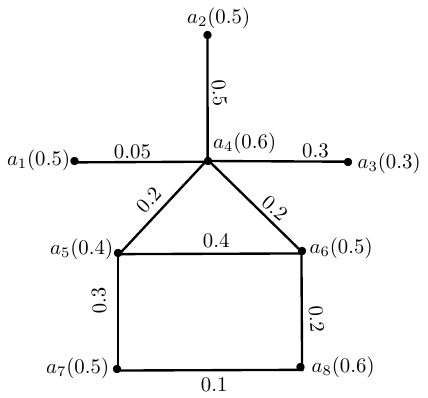} \caption{Illustration of SSD of a vertex in a FG.}
    \label{2}
\end{figure}
\begin{exam}\label{3}
For the FG in Figure \ref{2}, all edges are strong except $a_7a_8$. The MSDS containing $a_1$ with minimum weight is $D_1=\{a_1,a_2,a_3,a_5,a_6\}$ and $W(D_1)=1.35.$. Hence $_{sd}d(a_1)=1.25.$ Similarly, $_{sd}d(a_2)=_{sd}d(a_3)=1.35.$ And, $_{sd}d(a_4)=_{sd}d(a_5)=_{sd}d(a_6)=_{sd}d(a_8)=0.45$, $_{sd}d(a_7)=0.55.$ Therefore, $_{sp}\delta(G)=0.45$ and $_{sp}\Delta(G)=1.25.$ 
\end{exam}
\begin{rem}\label{4}
In a fuzzy star graph $FS_n$, all pendant vertices of the unique MST have SDD equal to $q$, the size of $FS_n$. And the vertex with degree $n$ has SDD equal to the minimum of weight of the edges of $FS_n$.\\Similarly for a CFG, each vertex has SDD equal to the minimum of weight of the vertices in the CFG.
\end{rem}
\begin{defn}\label{5}
A FG is $w-$strong domination regular fuzzy graph (w-SDRFG) or simply strong domination regular fuzzy graph (SDRFG), if $_{sd}d(u)=w$, $\forall u\in \mathcal{V}(\mathcal X)$. 
\end{defn}
\begin{exam}\label{6}
From Remark \ref{4}, for a CFG with n vertices $a_1,a_2,...,a_n$, $_{sp}d(a_i)=min \,\{\rho(a_j): 1\leq j\leq n\}$, $\forall 1\leq i\leq n.$ Hence, any CFG is a SDRFG.  
\end{exam}
Another example for a SDRFG is provided in Example \ref{8}.
\begin{figure}[h]
    \centering
    \includegraphics[height=4cm,width=11cm]{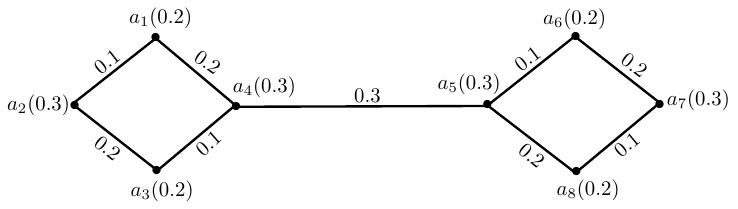}
    \caption{Strong domination regular fuzzy graph.}
    \label{7}
\end{figure}
\begin{exam}\label{8}
For the FG in Figure \ref{7}, SDD of every vertex is $3\times 0.1=0.3.$ Hence, the FG in Figure \ref{7} is a 0.3-SDRFG.
\end{exam}
Now, the concept of upper domination number is introduced using the weight of strong edges in Definition \ref{9}.
\begin{defn}\label{9}
 Let $\mathcal{X}$ be a FG. The upper strong domination number (USDN) of $\mathcal {X}$ is the maximum of weight of all the MSDSs of $\mathcal {X}$. The USDN of a FG is denoted by $\Gamma_s(\mathcal X)$, or simply $\Gamma_s.$
\end{defn}
An example to illustrate the concept of USDN is provided in Example \ref{11}.
\begin{figure}[h]
    \centering
    \includegraphics[height=5cm,width=5cm]{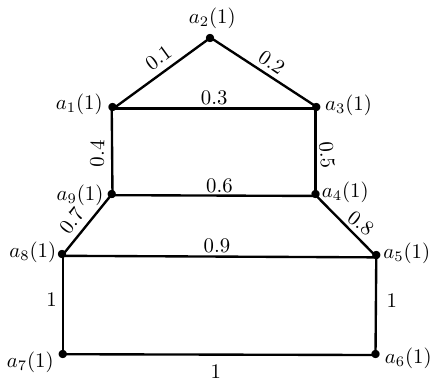}
    \caption{Illustration of USDN.}
    \label{10}
\end{figure}
\begin{exam}\label{11}
 For the FG in Figure \ref{10}, the SDS $D_1=\{a_3,a_6,a_9\}$ has the least weight. Hence, $\gamma_s= 0.2+1+0.4=1.6.$ Among the minimal SDSs the set $D_2=\{a_1,a_2,a_4,a_6,a_8\}$ has the maximum weight. Therefore, $\Gamma_s=0.2+0.5+1+0.4+0.7=2.8.$  
\end{exam}
Clearly, Theorem \ref{12} follow from the definitions of SDN and USDN.
\begin{thm}\label{12}
 For a FG $\mathcal{X}$, $\gamma_s(\mathcal{X})\leq \Gamma_s(\mathcal{X}).$
\end{thm}
\begin{thm}\label{13}
Let $\mathcal{X}$ be a FG, then 
$$\gamma_s(\mathcal{X})\leq _{sp}d(u)\leq \Gamma_s(\mathcal{X})$$ 
$\forall u \in V(\mathcal{X}).$
\end{thm}
\begin{proof}
The SDN is the minimum of weight of minimal SDSs of $\mathcal{X}$. Also, the USDN is the maximum of weight of the minimal SDSs of $\mathcal{X}$. The SDD of a vertex $u$ is the minimum of weight of minimal SDSs containing $u$. Hence from all the definitions the inequalities follow. 
\end{proof}
 The inequalities in Theorem \ref{13} is sharp for the SDRFGs in Figure \ref{7}.
The definition of irredundance number, upper irredundance number, independent domination number and independence number using weight of strong edges is defined as follows:
\begin{defn}\label{43} 
The minimum of weight of maximal irredundant sets of a FG, $\mathcal{X}$ is called strong irredundance number $ir_s(\mathcal{X})$.  
\end{defn}
\begin{defn}\label{44}
The maximum of weight of irredundant sets of a FG, $\mathcal{X}$ is called strong upper irredundance number $IR_s(\mathcal{X})$.
\end{defn}
\begin{defn}\label{45}
The minimum of weight of independent dominating sets in FG, $\mathcal{X}$ is called strong independent domination number $i_s(\mathcal{X})$.     
\end{defn}
\begin{defn}\label{46}
The maximum of weight of independent sets in FG, $\mathcal{X}$ is called strong independence number $\beta_s(\mathcal{X}).$  
\end{defn}
\begin{thm}\cite{nagoor2006}
 Every maximal strong independent set is a MSDS.   
\end{thm}
\begin{thm}\cite{nagoor2008}
A SDS in a fuzzy graph $\mathcal{X}$ is a MSDS iff it is strong dominating and fuzzy irredundant set.   
\end{thm}
\begin{thm}\cite{nagoor2008}\label{52}
Every MSDS is maximal fuzzy irredundant set.    
\end{thm}
\begin{thm}\cite{nagoor2011}\label{53}
Let $\mathcal{X}$ be a FG. If $x$ is a vertex not 
 strongly dominated by a maximal strong irredundant set $\mathscr{M}$. Then for some $a\in \mathscr{M}$, 
 \begin{enumerate}
     \item $P_{N_s}[a,\mathscr{M}]\subseteq N_s(x)$.
     \item For $a_1,a_2\in P_{N_s}[a,\mathscr{M}],$ such that $a_1\neq a_2$, either $b_i\in \mathscr{M}\setminus \{a\}$, for $i=1,2$ such that $a_i$ is strong neighbor of each vertex in $P_{N_s}[a_i,\mathscr{M}]$ or $a_1a_2$ is a strong edge.
 \end{enumerate}
\end{thm}
\begin{thm}\label{51}
Let $\mathcal{S}(\neq \phi) $ be the set of vertices in FG, $\mathcal{X}$ such that no vertex of $\mathcal{S}$ is strongly dominated by $\mathscr{M}$, a fuzzy maximal irredundant set with the least weight. Define $\mathscr{M}_x=\{a\in \mathscr{M}| P_N[a,\mathscr{M}]\subseteq N_s(x)\}$ for each $x\in \mathcal{S}.$ Let $\mathscr{B}\subseteq \mathscr{M}$ with the least weight $w$ such that $\mathscr{M}_x\cap \mathscr{B} \neq \phi $ for each $X\in \mathcal{X}$. Then $\gamma_s< ir_s+w.$   
\end{thm}
\begin{proof}
 Let $\mathscr{M}=\{a_1,a_2,...,a_n\}$ be a strong maximal irredundant set with the least weight. Therefore, $ir_s=W(\mathscr{M}).$ Let $\mathscr{B}=\{b_1,b_2,...,b_m\}$, $n\leq m.$ BY definition of $\mathscr{B}$, for each $j=1,2,...,m$, $b_j\in \mathscr{M}_x$ for some $x\in S.$ Hence $P_{N_s}[b_j,\mathscr{M}]\subseteq N(x).$ Since no vertex of $\mathscr{M}$ strongly dominates $x,$ $b_j\notin N(x)$, which implies $b_j\notin P_{N_s}[b_j,\mathscr{M}].$ Hence $\exists$ $b_j'\in P_{N_s}[b_j,\mathscr{M}], \, b_j'\neq b_j .$ Therefore, $\exists \, b_j', 1\leq j\leq m $ such that $b_j'\in N_s(x) $for each $x\in S.$  Now consider the set $\mathcal{D}=\mathscr \cup \{b_1',b_2',...,b_m'\}$. Here, $\mathcal{D}$ is a SDS with weight $ir_s+w$. Since $\mathscr{M}\subset \mathcal{D}$, by maximality of $\mathscr{M}$, $\mathcal{D}$ is not irredundant. Therefore, by Theorem \ref{52}, $\mathcal{D}$ properly contains a minimal SDS. Hence, $\gamma< ir_s+w.$
\end{proof}
From definition of SDN and definitions \ref{1}, \ref{9}, \ref{43}, \ref{44}. \ref{45}, and \ref{46} the inequality in Theorem \ref{47} follows:
\begin{thm}\label{47}
Let $\mathcal{X}$ be a FG, then $$ir_s\leq \gamma_s \leq i_s\leq \beta_s \leq \Gamma_s \leq IR_s.$$ 
\end{thm}
Now, an example to illustrate definitions \ref{43}, \ref{44}, \ref{45} and \ref{46}, \ref{9} and $\gamma_s$ is provided in Example \ref{49}.
\begin{figure}[h]
    \centering
    \includegraphics[height=5cm,width=5cm]{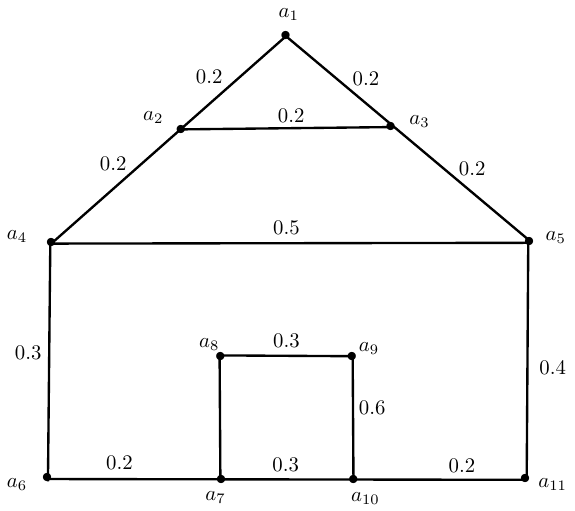}
    \caption{Illustration for Theorem \ref{47}}
    \label{48}
\end{figure}
\begin{exam}\label{49}
 Consider FG in Figure \ref{48}. Here $\{a_2,a_3,a_7,a_{10}\}$ is a maximal irredunant set with least weight. Hence $ir_s=0.8.$ The set $\{a_2,a_4,a_7,a_{10}\}$ is a minimum SDS, hence $\gamma=0.8.$ The set $\{a_1,a_4,a_{11},a{7}\}$ is an independent set with minimum weight the least weight. Therefore, $i_s=0.9$. Also an independent set with maximum weight is $\{a_2,a_5,a_8,a_{10},a_6\}$, thus $\beta_s=1.3.$ A MSDS with maximum weight is $\{a_1,a_4,a_5,a_8,a_9\}$ which is also an irredundant set with maximum weight. Hence $\Gamma_s=IR_s=1.4.$ 
\end{exam}
From \ref{13} and \ref{47} the Proposition \ref{prop} follows:
\begin{prop}\label{prop}
Let $\mathcal{X}$ be a FG, then 
$$ir_s\leq \gamma \leq _{sp}d(u)\leq \Gamma \leq IR_s$$ 
$\forall u \in V(\mathcal{X}).$    
\end{prop}
\section{Strong Domination Index in Fuzzy Graphs}
In this section the notion of Strong domination index is defined in FGs. Inequalities for SDI of a FG is obtained. The idea of SDI is discussed in various FGs like CFGs, CBFGs, complete $r-$ partite FGs, fuzzy trees, fuzzy cycles and fuzzy stars. The union and join of FGs are also considered in the section.
\begin{defn}\label{14}
Let $\mathcal{X}$ be a FG. Then the strong domination index (SDI) of $\mathcal{X}$, denoted as $SDI(\mathcal{X})$ is defined as the sum of SDD of vertices of $\mathcal{X},$ i.e,
$$SDI(\mathcal{X})=\sum\limits_{v\in \mathcal{V}(\mathcal{X})} {_{sd}d(v)}.$$
\end{defn}
\begin{exam}\label{15}
For the FG in Figure \ref{2}, the SDI is 
$SDI(\mathcal{X})={_{sd}d(a_1)}+{_{sd}d(a_2)}+{_{sd}d(a_3)}+...+{_{sd}d(a_8)}=3\times1.25+4\times0.45+0.55
=6.1$
\end{exam}
From Theorem \ref{12}, $\gamma(\mathcal{X})\leq _{sp}d(u)\leq \Gamma(\mathcal{X})$ $\forall u \in V(\mathcal{X}).$ Hence, $\sum\limits_{u\in V(\mathcal{X})}\gamma(\mathcal{X})\leq \sum\limits_{u\in V(\mathcal{X})} {_{sp}d(u)}\leq \sum\limits_{u\in V(\mathcal{X})}\Gamma(\mathcal{X})$. Therefore, $n\gamma(\mathcal{X})\leq SDI(\mathcal{X})\leq n\Gamma(\mathcal{X})$. And hence the theorem follows:\\
\begin{thm}\label{16}
    For a FG $\mathcal{X}$, with SDN $\gamma$ and USDN $\Gamma$, $$\gamma(\mathcal{X})\leq \frac{SDI(\mathcal{X})}{n}\leq \Gamma(\mathcal{X})$$,
    where $n=|\rho^*|.$
\end{thm}
From Proposition \ref{prop}, \ref{prop2} follows:
\begin{prop}\label{prop2}
 Let $\mathcal{X}$ be a FG with $|\mathcal{V}(\mathcal{X})|=n$, then 
$$ir_s\leq \gamma \leq \frac{SDI(\mathcal{X})}{n}\leq \Gamma \leq IR_s.$$    
\end{prop}
\begin{thm}\label{17}
Let $\mathcal{X}$ be a CFG with $n$ vertices. Then, $$SDI(\mathcal{X})= n \rho(u),$$
where $u$ is the vertex with least weight, and $n=|\rho^*|.$.
\end{thm}
\begin{proof}
For a CFG $\mathcal{X}$, each set $\{v\},$ $v\in V(\mathcal{X})$ is a SDS. And, the weight of each SDS ${v}$ is $\rho(u)$, where $u$ is the vertex with the least weight. Hence the SDD $_{sd}d(v)$ is $\rho(u)$, $\forall v\in V(\mathcal{X}).$ Therefore, $SDI(\mathcal{X})=n \rho(u) .$
\end{proof}
\begin{rem}\label{25}
In all cases of complete bipartite FG and complete $r-$ partite FGs, take $N_1$ as the partite set containing the least weight vertex.    
\end{rem}
\begin{thm}\label{18} 
 Let $\mathcal{X}$ be a complete bipartite FG with partitions $N_1$ and $N_2$, $|N_1|\geq 1, |N_2|\geq 1$. Let the vertices in each partition be labelled as $a_1,a_2,...,a_{n_1}$ and $b_1,b_2,...,b_{n_2}$ respectively. Also, let the $k^{th}$ minimum weight $b_i\in N_2$, then
 $$SDI(\mathcal{X})=[2(n_2+1)+(n_1-1)]\rho(a_1)+\rho(a_2)+...+\rho(a_{k-1})+(n_1-(k-1)\rho(b_i)$$ where $a_1$ is the vertex with least weight. 
\end{thm}
\begin{proof}
  Let $\mathcal{X}$ be a complete bipartite FG with partitions $N_1$ and $N_2$, $|N_1|\geq 1, |N_2|\geq 1$. Since the least weight vertex belongs to $N_1$, for all vertices $b_j, j=1,2,...,n_2$ in $N_2$, the minimum of weight of strong edges incident at $b_j$ is $\rho(a_1)$. Hence, the contribution of vertex $b_j$ in $N_2$ to the weight of SDS containing $b_j$ is $\rho(a_1)$. Similarly, since the $k^{th}$ minimum weight vertex belongs to $N_2$, the minimum of weight of edge incident at vertices $a_2,a_3,...,a_{k-1}$ is $\rho(a_2),\rho(a_3),...,\rho(a_{k-1})$ respectively. Hence, the contribution of vertices $a_2,a_3,...,a_{k-1}$ to the weight of SDS containing $a_2,a_3,...,a_{k-1}$ is $\rho(a_2),\rho(a_3),...,\rho(a_{k-1})$ respectively. And, for all other vertices $a_k,a_{k+1},...,a_{n_1}$, the minimum of weight of edges incident at $a_k,a_{k+1},...,a_{n_1}$ is $\rho(b_i)$. For any vertex $b_j$ in $N_2$, the MSDS containing $b_j$ is $\{b_j, a_1\}$ with weight $2\rho(a_1)$, i.e, SDD of $b_j$ is $2\rho(a_1)$ . Similarly, for any vertex $a_l$ in $N_1$, a MSDS containing $a_l$ is $\{a_l,b_1\}$. For vertex $a_1$, SDD is $2\rho(a_1)$. And for vertices $a_2,a_3,...,a_{k-1}$, SDD is $\rho(a_2)+\rho(a_1),\rho(a_3)+\rho(a_1),...,\rho(a_{k-1})+\rho(a_1)$. The SDD of the remaining $(n_1-(k-1))$ vertices in $N_1$ is $\rho(b_i)+\rho(a_1)$. Therefore, the SDI of $\mathcal{X}$ is $[2(n_2+1)+(n_1-1)]\rho(a_1)+\rho(a_2)+...+\rho(a_{k-1}+(n_1-(k-1)\rho(b_i)$.  
\end{proof}
\begin{defn}\label{19}
  A FG $\mathcal{X}=(\rho,\upsilon)$ is said to be $r-$partite if the vertex set $\mathcal{V}$ can be partitioned into $r$ non-empty sets $\mathcal{V}_1, \mathcal{V}_2,..., \mathcal{V}_r$ such that $\upsilon(a,b)=0$ if $a,b\in \mathcal{V}_i$, for $i=1,2,...,r$. Also, if $\upsilon(a,b)=\rho(a)\wedge\rho(b)$ for all $a\in \mathcal{V}_i$ and $b\in V_j, i\neq j$, then $\mathcal{X}$ is called a complete $r-$ partite FG and is denoted by $K_{\rho_1,\rho_2,...,\rho_r}.$  
\end{defn}
\begin{cor}\label{20}
Let $\mathcal{X}$ be a complete r-partite FG with partitions $N_1$,$N_2$,...,$N_k$, $|N_1|\geq 1, |N_2|\geq 1$,..., $|N_k|\geq 1$ . Let the vertices in each partition be labelled as $a_{1_1},a_{1_2},...,a_{1_{n_1}}$, $a_{2_1},a_{2_2},...,a_{2_{n_2}}$,..., $a_{k_1},a_{k_2},...,a_{k_{n_k}}$ respectively. Also, let the $k^{th}$ minimum weight $b_i\in N_m, m\neq 1$, then
 $$SDI(\mathcal{X})=[2(n_2+n_3+...+n_k+1)+(n_1-1)]\rho(a_1)+\rho(a_2)+...+\rho(a_{k-1})+(n_1-(k-1)\rho(b_i)$$ where $a_1$ is the vertex with least weight.
\end{cor}
\textbf{Illustration of SDI in complete bipartite FGs and complete $r-$ partite FGs}
\vspace{0.5cm}

\begin{figure}[!h]
\begin{minipage}[c]{0.4\linewidth}
    \centering
     \vspace{2cm}
     
    \includegraphics[height=5cm,width=4cm]{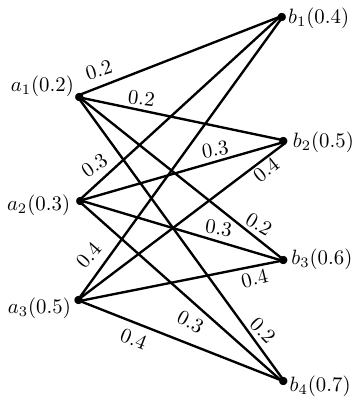}
    \caption{Complete bipartite fuzzy graph}
    \label{22}
\end{minipage}\hfill
\begin{minipage}[c]{0.4\linewidth}
    \centering
    \includegraphics[height=7cm,width=4cm]{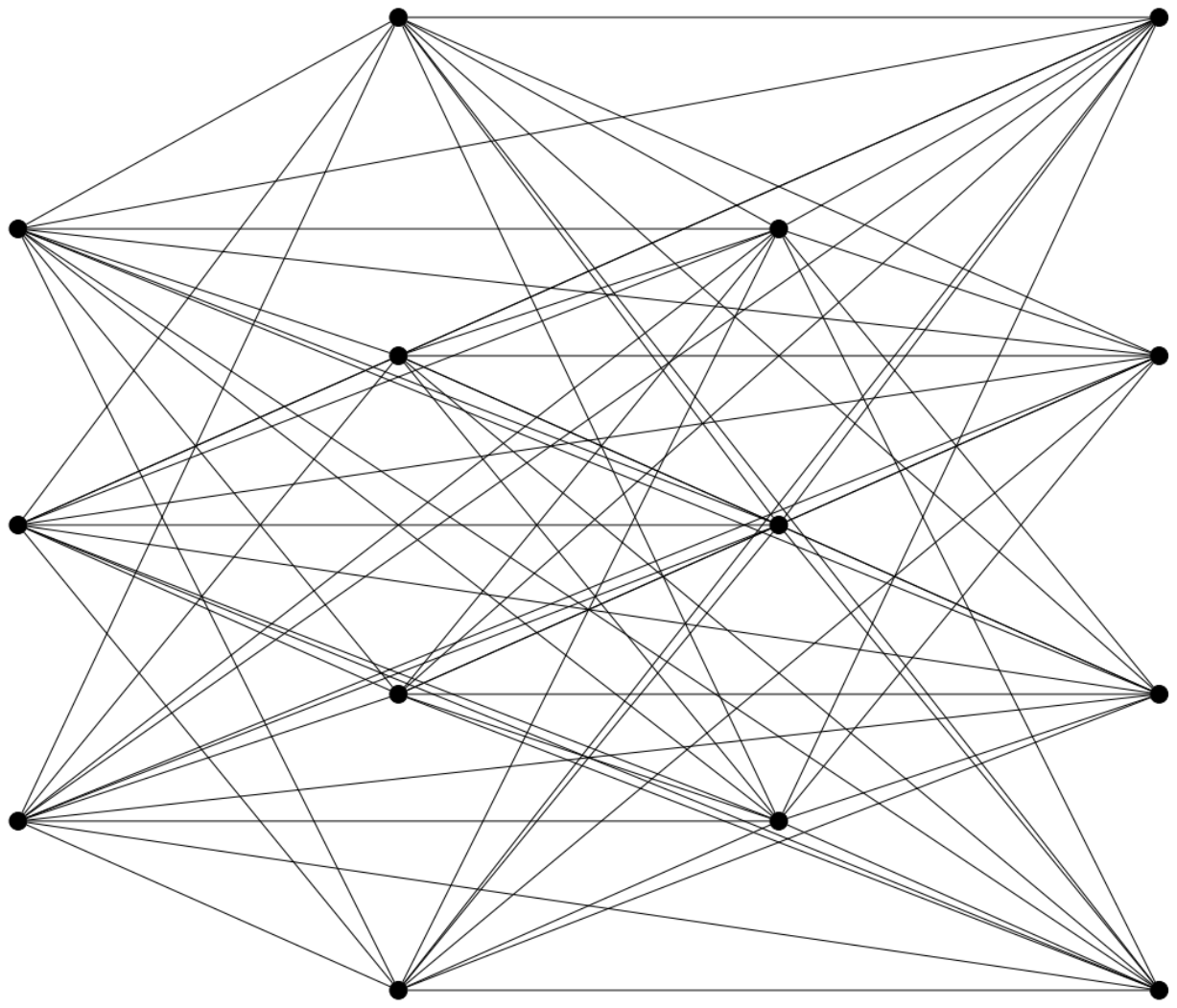}
    \caption{Complete 4-partite graph}
    \label{23}
\end{minipage}
\end{figure}

\begin{exam}\label{21}
Consider the complete bipartite fuzzy graph in Figure \ref{22}. The partitions are $N_1$ and $N_2$. $|N_1|=3$, $|N_2|=4$ and $N_1=\{a_1,a_2,a_3\}$, $N_2=\{b_1,b_2,b_3,b_4\}$. The minimum weight vertex is $a_1$ with weight 0.2. The $3^{rd}$ minimum weight vertex $b_1$ belongs to $N_2$ with weight 0.4. The SDDs of $a
_1$, $a_2$ and $a_3$ are $0.2+0.2$, $0.3+0.2$ and $0.4+0.2$ respectively with the MSDS contributing to the SDD being $\{a_1,b_1\}$, $\{a_2,b_1\}$ and $\{a_3,b_1\}$ respectively. Similarly, the SDDs of vertices $b_1,b_2,b_3,b_4$ is $0.2+0.2$ and the MSDSs contributing to the SDD are $\{b_1,a_1\},\{b_2,a_1\},\{b_3,a_1\}$ and $\{b_4,a_1\}$
Hence, the SDI of FG in Figure \ref{23} is $0.2+0.2+0.3+0.2+0.4+0.2+4(0.2+0.2)$, i.e, $[2(4+1)+(3-1)]0.2+0.3+0.4= [2(n_2+1)+(n_1-1)]0.2+0.3+(n_1-(k-1))0.4.$
\end{exam}
\begin{exam}\label{24}
 Consider a FG with the underlying graph as in Figure \ref{23}. The vertex set $V$ is partitioned into 4 sets $N_1,N_2,N_3,$ and $N_4$ with $|N_1|=3,|N_2|=4, |N_3|=3$ and $|N_4|=4$. Let $N_1=\{a_{1_1},a_{1_2},a_{1_3}\}$,  $N_2=\{a_{2_1},a_{2_2},a_{2_3},a_{2_4}\}$, $N_3=\{a_{3_1},a_{3_2},a_{3_3}\}$ and $N_4=\{a_{4_1},a_{4_2},a_{4_3},a_{4_4}\}$. The weights are given by $\rho(a_{1_1})=0.1, \rho(a_{1_2})=0.2, \rho(a_{1_3})=0.5, \rho(a_{2_1})=0.4, \rho(a_{2_2})=0.5, \rho(a_{2_3})=0.6, \rho(a_{2_4})=0.4, \rho(a_{3_1})=0.7, \rho(a_{3_2})=0.5, \rho(a_{3_3})=0.5, \rho(a_{4_1})=0.3, \rho(a_{4_2})=0.5, \rho(a_{4_3})=0.6$ and $\rho(a_{4_4})=0.7.$ For a vertex $a$ in $N_2,N_3,N_4$ the SDD is $0.1+0.1$ and the MSDS contributing to the SDD is $\{a,a_{1_1}\}$. Similarly, for vertex $a'$ a MSDS contributing to the SDD is $\{a',a_{4_1}\}$. The SDD is $\rho(a')+0.1$ if $\rho(a')\leq\rho(a_{4_1})$ and $\rho(a_{4_1})+0.1$ if $\rho(a')\geq\rho(a_{4_1})$. Hence, SDI of the FG is $[2(4+3+4+1)+2]0.1+0.2+0.3=3.1$, i.e, $[2(n_2+n_3+n_4+1)+(n_1-1)]\rho(a_{1_1})+\rho(a_{1_2})+(n_1-(k-1))\rho(a_{4_1}).$
\end{exam}
\begin{thm}\label{26}
    The SDI of a fuzzy star graph $FS_n$ is
    $$SDI(FS_n)=\wedge \{\upsilon(uv): uv\in \upsilon^*\}+(n)q$$
\end{thm}
\begin{proof}
    Let $FS_n$ be the fuzzy star graph with $n+1$ vertices $a_1,a_2,...,a_{n+1}.$ Let $a_1$ be the vertex with $|N_s|=n$. Then the SDD of $a_1$ is the minimum of weight of the edges and the MSDS contributing to the SDD is $\{a_1\}$. Now, for all remaining vertices $a_2,a_3,...a_{n+1}$, the MSDS contributing to the SDD is the set $\{a_2,a_3,...,a_{n+1}\}$ with SDD equal to the sum of weight of all the edges, i.e, $q.$ Hence, $SDI(FS_n)=\wedge \{\upsilon(uv): uv\in \upsilon^*\}+(n)q$. 
\end{proof}
\begin{prop}\label{27}
In a fuzzy cycle all $\beta-$ edges have equal weight.    
\end{prop}
\begin{prop}\label{28}
In a fuzzy cycle the weight of a $\beta-$ edge is always less than or equal to the weight of an $\alpha-$ edge with the least weight
\end{prop}
\begin{proof}
 Let $\mathcal{X}$ be a fuzzy cycle. Let $ab$ be an $\alpha-$ edge with least weight. Suppose $uv$ is a $beta edge$, and $\upsilon(ab)<\upsilon(uv).$ Then, $CONN_{\mathcal{X}\ uv}(u,v)=\upsilon(ab)$, which is contradiction to the fact that $\mathcal{X}$ is a fuzzy cycle. Hence the weight of a $\beta-$ edge is always less than or equal to the weight of an $\alpha-$ edge with the least weight.
\end{proof}
\begin{thm}\label{29}
    If the fuzzy cycle $\mathcal{X}$, is $\beta$ saturated and the weight of the $\beta-$ edge is $w$, then $$SDI(\mathcal{X})= n \bigg\lceil{\frac{n}{3}}\bigg\rceil w$$
\end{thm}
\begin{proof}
    Let $\mathcal{X}$ be a fuzzy cycle. Since $\mathcal{X}$ is $\beta$ saturated there is a $\beta-$ edge incident at every vertex. Hence for every vertex $a$ the minimum of weight of edges incident at $a$ is the weight of $\beta-$ edge. By \ref{27} each vertex contributes the same weight say $w$ to the weight of SDS. Also the domination number of a cycle is $\bigg\lceil{\frac{n}{3}}\bigg\rceil$. Hence, the SDD of each vertex is $\bigg\lceil{\frac{n}{3}}\bigg\rceil w$. Hence, $SDI(\mathcal{X})=n \bigg\lceil{\frac{n}{3}}\bigg\rceil w.$
\end{proof}
\begin{thm}\label{30}
    Let $\mathcal{X}$ be a FG. If $ab$ is a $\delta-$ edge, then $SDI(\mathcal{X})= SDI(\mathcal{X}\setminus ab).$
\end{thm}
\begin{proof}
    By the definition of strong domination it is clear that a vertex $a$ dominates a vertex $b$  if $ab$ is strong edge. And also only the weight of strong edges is considered while calculating SDN. Hence, the deletion of a $\delta-$ has no effect on SDN of the FG $\mathcal{X}.$ Therefore, if $ab $ is a $\delta-$ edge, $SDI(\mathcal{X})= SDI(\mathcal{X}\setminus ab).$ 
\end{proof}
\begin{thm}\label{32}
 Let $\mathcal{X}$ be a fuzzy tree and ${\mathcal{X}}^*$ is not a tree. Then there exist at least one edge $uv\in \upsilon^*$ such that $SDI(\mathcal{X}\setminus uv)=SDI(\mathcal{X}).$
\end{thm}
\begin{proof}
Let $\mathcal{X}$ be a fuzzy tree and ${\mathcal{X}}^*$ is not a tree. Then there exists a cycle in $\mathcal{X}$ and an edge $uv$ in the cycle such that $\upsilon(uv)<CONN_{\mathcal{X}\setminus uv}(u,v)$, i.e, $uv$ is a $\delta-$ edge. Since the deletion of  $\delta-$ edges have no effect on the strong domination $SDI(\mathcal{X}\setminus uv)=SDI(\mathcal{X}).$ Let $\mathcal{Y}=\mathcal{X}\setminus uv$, and if $\mathcal{Y}^*$ is not a tree, then there exists more edges like $uv$ in $\mathcal{X}$. Hence there exist at least one edge $uv\in \upsilon^*$ such that $SDI(\mathcal{X}\setminus uv)=SDI(\mathcal{X}).$  
\end{proof}
\begin{thm}\label{31}
    Let $\mathcal{X}$ be a fuzzy tree, if $\mathcal{F}$ is the unique maximum spanning tree of $\mathcal{X}$, then $SDI(\mathcal{X})=SDI(\mathcal{F}).$
\end{thm}
\begin{proof}
     The proof of Theorem 2.3.1 \cite{book2} discusses the construction of $\mathcal{F}$. $\mathcal{F}$ is obtained from $\mathcal{X}$ by deleting the edge $xy$ of a cycle if $\upsilon(x,y)<CONN_{\mathcal{X}\setminus xy }(x,y)$, i.e, by deleting all the $\delta-$ edges of $\mathcal{X}$. By \ref{30}, the deletion a $\delta-$ edge has no effect on the SDI of $\mathcal{X}$. Applying \ref{30} to each $\delta-$ edge in $\mathcal{X}$, the SDI of $\mathcal{F}$ is obtained as $SDI(\mathcal{X})=SDI(\mathcal{F}).$
\end{proof}
\begin{prop}\label{42}
Let $\mathcal{X}$ be a connected strong FG such that $\rho(u)=1\, \forall u\in \rho* $ and  $|\rho^*|>1$, then $SDI(\mathcal{X})\leq |\rho*| WI(\mathcal{X}).$
\end{prop}
\begin{proof}
  Let $\mathcal{X}$ be a connected strong FG such that $\rho(u)=1\, \forall u\in \rho* $. Consider vertex $u$ and a MSDS $D_u$ containing $u$ which is considered for $_{sd}d(u)$. \\
  \textbf{Case 1:} $D_u=\{u\}$.\\
  Since $\mathcal{X}$ is connected, $\exists $ at least one vertex $v$ such that $uv$ is an edge in $\mathcal {X}$. Hence, $_{sd}d(u)\leq WI(\mathcal{X}).$\\
  \textbf{Case 2:} $D_u=\{u,u_1\}$.\\
  In this case either $uu_1\in \upsilon^*$ or $uu_1\notin \upsilon^*$
  If $uu_1\in \upsilon^*$, then since $\{u,u_1\}$ is MSDS and $\mathcal{X}$ is connected $uu_1$ is not the only edge in $\mathcal{X}$ and there exists at least two distinct vertices that are private neighbors of $u$ and $u_1$ respectively.  If $uu_1\notin \upsilon^*$, since $\mathcal{X}$ is connected $\exists$ at least one vertex $v$ in path from $u$ to $u_1$. Hence, $\exists$ at least two edges in the path. Also, each vertex $a$ contributes the minimum of weight of edges incident at $a$ to the weight of $D_u$. Hence $_{sd}d(u)\leq WI(\mathcal{X}).$\\
  \textbf{Case 3:} $D_u=\{u,u_1,u_2,...,u_n\}$.\\
  For each vertex in $a\in\{u,u_1,u_2,...,u_n\}$ minimum of weight of edges incident at $a$ is contributed towards weight of the MSDS $D_u=\{u,u_1,u_2,...,u_n\}.$ Each vertex $a$ has either a private neighbor $b$ such that $b\neq a$ or $P_N[a,D_u]=\{a\}$, which means $a $ is isolated in $\mathcal{X}[D_u]$. Among $u,u_1,u_2,...,u_n$ let $a_1,a_2,...,a_k$ be the vertices with at least one  distinct vertex in the private neighborhood, i.e, $\exists$ vertices $b_1,b_2,...,b_k$ such that $b_i\in P_N[a_i,D_u], 1\leq i\leq k .$ And let $c_1,c_2,...,c_{k'}$, $k+k'=n+1$ be the vertices such that $P_N[c_i,D_u]=\{c_i\}, 1\leq i\leq k'.$ For vertex $a_i$ the minimum of weight of edges incident at $a_i$ is less than $\upsilon(a_i b_i), 1\leq i\leq k.$ Now, consider two vertices among $c_1,c_2,...,c_{k'}$. A path from $c_i$ to $c_j$, $1\leq i,j\leq k', i\neq j$ contains at least one vertex $v_{i'}$ such that $v_{i'}\notin N_s[a],$ $a\in \{b_1,b_2,...,b_k\} $. Therefore, $\exists$ at least two edges in each path. 
  Suppose there is only one vertex $b$ such that $P_N[b, D_u]=\{b\}$. Then $b$ is isolated in $\mathcal{X}[D_u]$. Since $\mathcal{X}$ is connected, consider any path from $b$ to $u_i$, then $\exists$ at least two edges in such paths. And the edge weight contributed by $b$ is distinct from $\upsilon(a_ib_i).$ Hence $_{sd}d(u)\leq WI(\mathcal{X}).$
\end{proof}
\begin{thm}\label{33}
    Let $\mathcal{X}_1=(\rho_1,\upsilon_1)$ and $\mathcal{X}_2=(\rho_2,\upsilon_2)$ be isomorphic. Then $SDI(\mathcal{X}_1)=SDI(\mathcal{X}_2).$
\end{thm}
\begin{proof}
    Let $\mathcal{X}_1=(\rho_1,\upsilon_1)$ and $\mathcal{X}_2=(\rho_2,\upsilon_2)$ be isomorphic and $\theta$ be the bijection from ${\rho_1}^*$ to ${\rho_2}^*$ such that $\rho_1(u)=\rho_2(\theta(u))$ for $u\in{\rho_1}^*$ and $\upsilon(u,v)=\upsilon(\theta(u),\theta(v))$ for $uv\in {\upsilon_1}^*$. Since $\mathcal{X}_1$ and $\mathcal{X}_2$ are isomorphic the weight of edges incident at $u$ and $\Theta(u)$ are same. Hence, the minimum of weight of edges incident at $u$ and $\theta(u)$ are same. Hence $_{sd}d(u)=_{sd}d(\theta(u))$. 
    \begin{align*}
     SDI(\mathcal{X}_1)&=\sum\limits_{v\in V(\mathcal{X}_1)} {_{sd}d(v)} \\
     &= \sum\limits_{\theta(v)\in V(\mathcal{X}2)} {_{sd}d(\theta(v))}\\
     &=SDI(\mathcal{X}_2)
    \end{align*}
     Therefore, $SDI(\mathcal{X}_1)=SDI(\mathcal{X}_2).$
\end{proof}
\begin{thm}\label{34}
    Let $\mathcal{X}=\bigcup\limits_{i=1}^t \mathcal{X}_i$ be the disjoint union of FGs $\mathcal{X}_1, \mathcal{X}_2,...,\mathcal{X}_t$, then $$SDI(\mathcal{X})={\sum\limits_{i=1}^t SDI(\mathcal{X}_i)} + \sum\limits_{i=1}^t \big(n\setminus |\mathcal{V}(\mathcal{X}_i)|\big)\gamma(\mathcal{X}_i),$$ where $n=|\mathcal{V}(\mathcal{X})|$ and $|\mathcal{V}(\mathcal{X}_i)|$ is the number vertices in each component $\mathcal{X}_i.$
\end{thm}
\begin{proof}
    Let $\mathcal{X}=\bigcup\limits_{i=1}^t \mathcal{X}_i$ be the disjoint union of FGs $\mathcal{X}_1, \mathcal{X}_2,...,\mathcal{X}_t$. Consider any component $\mathcal{X}_k$ of  $\mathcal{X}$ and let $u$ be a vertex of $\mathcal{X}_k$. A SDS of $\mathcal{X}$ is disjoint union of SDSs of each component $\mathcal{X}_i$. Let $D$ be the MSDS considered for the $_{sd}d(u)$.  For all other components other than $\mathcal{X}_k$, the SDS contributed from each $\mathcal{X}_k$ is the minimum SDS, since it is the least weight SDS satisfying the minimality condition. And, the SDS contributed from $\mathcal{X}_k$ is the SDS considered for $_{sd}d(u)$ in $\mathcal{X}_k$. Hence, by adding $_{sd}d(u)$ of each vertex in $\mathcal{X}_k$, the $SDI(\mathcal{X}_k)$ is obtained with the dominating number of each other components $\mathcal{X}_k$ being added $n\setminus |V(\mathcal{X}_k)|$ times. Therefore $$SDI(\mathcal{X})={\sum\limits_{i=1}^t SDI(\mathcal{X}_i)} + \sum\limits_{i=1}^t \big(n\setminus |\mathcal{V}(\mathcal{X}_i)|\big)\gamma(\mathcal{X}_i).$$
\end{proof}
\begin{thm}\label{35}
Let $\mathcal{X}_1$ and $\mathcal{X}_2$ be two FGs with $m$ and $n$ vertices respectively. Let $\mathcal{X}=\mathcal{X}_1+\mathcal{X}_2$ be the join of $\mathcal{X}_1$ and $\mathcal{X}_2$. If $D\subset V(\mathcal{X}_1)$ is not a SDS of $\mathcal{X}_1$, then $D$ is not a SDS of $\mathcal{X}$.
\end{thm}
\begin{proof}
    By definition of join of FGs the edges of the form $uv,$ where $u\in V(\mathcal{X}_1)$ and $v\in V(\mathcal{X}_2)$ are strong edges. Suppose an edge be of the form $u'v'$, where $u',v'\in \mathcal{X}_1$ or $u',v'\in \mathcal{X}_2$ is not strong. Assume that $u',v'\in \mathcal{X}_1$, then $\upsilon_1(u',v')<CONN_{\mathcal{X}_1\setminus u'v'}(u',v')$. Hence there exists a path, say $P$ in $\mathcal{X}_1$ with strength of connectedness greater than $\upsilon_1(u',v')$. By the definition of join of FGs, the path $P$ exists in $\mathcal{X}$. And since the strength of connectedness between $u'$ and $v'$ is the maximum of strength of connectedness of all paths between $u'$ and $v'$, $\upsilon(u',v')< CONN_{\tilde{\mathcal X}\setminus u'v'}(u',v')$. Therefore, $u'v'$ is not strong in $\mathcal{X}.$ By similar way it can be proved in the case when $u',v'\in \mathcal{X}_2$. Hence, if any edge $u'v'$, $u',v'\in \mathcal{X}_1$ or $u',v'\in \mathcal{X}_2$ is not strong, then $u'v'$ is not strong in $\mathcal{X}.$ So, if $D $ is not a SDS in $\mathcal{X}_i,\, i=1,2$, then $D$ is not a SDS $\mathcal{X}.$  
\end{proof}
\begin{thm}\label{36}
    Let $\mathcal{X}_1$ and $\mathcal{X}_2$ be two FGs with $m$ and $n$ vertices respectively. Let $\mathcal{X}=\mathcal{X}_1+\mathcal{X}_2$ be the join of $\mathcal{X}_1$ and $\mathcal{X}_2$. Then, 
    $$SDI(\mathcal{X})= \sum\limits_{\substack{u_i\in V(\mathcal{X}_1)\\ 1 \leq i \leq m}} \min\{W(\{u_i,v\}), W(\{D_{u_i})\}\} + \sum\limits_{\substack{v_i\in V(\mathcal{X}_2)\\ 1 \leq i \leq n}} \min\{W(\{u,v_i\}), W(\{D_{v_i})\}\}$$
    where $v$ in the first summation is vertex in $\mathcal{X}_2$ having the least weight and $u$ in the second summation is vertex in $\mathcal{X}_1$ having the least weight. Also, $D_{u_i}$ is any MSDS containing $u_i$ in $\mathcal{X}_1$ which is also a MSDS in $\mathcal{X}$ and $D_{v_i}$ is any MSDS containing $v_i$ in $\mathcal{X}_2$ which is also a MSDS in $\mathcal{X}$.
\end{thm}
\begin{proof}
Let $\mathcal{X}_1$ and $\mathcal{X}_2$ be two FGs with $m$ and $n$ vertices respectively. Let $\mathcal{X}=\mathcal{X}_1+\mathcal{X}_2$ be the join of $\mathcal{X}_1$ and $\mathcal{X}_2$. Consider the vertices $\{u_1,u_2,...,u_m\}$ of $\mathcal{X}_1$. Let $D_{u_i}, 1\leq i\leq m$ be a MSDS in $\mathcal{X}_1$ containing $u_i$ such that the edges connecting the vertices of $D_{u_i}$ to the vertices of $V(\mathcal{X})\setminus D_{u_i} $ remains strong in $\mathcal{X}.$ Hence $D_{u_i}$ dominates all vertices of $\mathcal{X}_1$ and any vertex of $\mathcal{X}_1$ dominates all vertices of $\mathcal{X}_2 $ in the join. Therefore, $D_{u_i}$ is a SDS in $\mathcal{X}$ also. Also by Theorem.\ref{35}, any subset $S\subset V(\mathcal{X}_1)$ which is not a SDS of $\mathcal{X}_1$ cannot be a SDS in $\mathcal{X}.$ Now, consider $\{u_i, v'\}$, where $v'$ is a vertex in $\mathcal{X}_2$. Here $u_i$ dominates all vertices of $\mathcal{X}_2$ and $v'$ dominates all vertices of $\mathcal{X}_1$ in $\mathcal{X}.$ Hence, $\{u_i, v'\}$ is also a SDS in $\mathcal{X}$ containing $u_i.$ Among all such sets the one with least weight is $\{u_i,v\}$, where $v$ is the least weight vertex in $\mathcal{X}_2.$ So, any MSDS in $\mathcal{X}$ containing $u_i$ is a MSDS in $\mathcal{X}_1$ containing $u_i$ or a SDS of the form $\{u_1,v\}, $ where $v\in V(\mathcal{X}_2)$. Therefore, $_{sd}d(u_i)=\sum\min\{W(\{u_i,v\}), W(\{D_{u_i})\}\}$ where $v$ is the least weight vertex in $\mathcal{X}_2.$ Similarly the case of vertices of $\mathcal{X}_2$ can be proved and $_{sd}d(v_i)=\sum\min\{W(\{u,v_i\}), W(\{D_{v_i})\}\}$ where $u$ is the least weight vertex in $\mathcal{X}_1.$ Hence, the $SDI(\mathcal{X})$ is, $$SDI(\mathcal{X})= \sum\limits_{\substack{u_i\in V(\mathcal{X}_1)\\ 1 \leq i \leq m}} \min\{W(\{u_i,v\}), W(\{D_{u_i})\}\} + \sum\limits_{\substack{v_i\in V(\mathcal{X}_2)\\ 1 \leq i \leq n}} \min\{W(\{u,v_i\}), W(\{D_{v_i})\}\}$$.
\end{proof}
\begin{thm}
 Let $\mathcal{X}_1$ and $\mathcal{X}_2$ be two strong FGs with $m$ and $n$ vertices respectively. Let $\mathcal{X}=\mathcal{X}_1+\mathcal{X}_2$ be the strong join of $\mathcal{X}_1$ and $\mathcal{X}_2$. Then, $SDI(\mathcal{X})\leq SDI(\mathcal{X}_1)+SDI(\mathcal{X}_2).$  
\end{thm}
\begin{proof}
  Let $\mathcal{X}_1$ and $\mathcal{X}_2$ be two strong FGs with $m$ and $n$ vertices respectively. Let $\mathcal{X}=\mathcal{X}_1+\mathcal{X}_2$ be the strong join of $\mathcal{X}_1$ and $\mathcal{X}_2$. Since the join is strong, every MSDS of $\mathcal{X}_1$ or $\mathcal{X}_2$ is a SDS of $\mathcal{X}$. Hence the MSDS containing vertex $a, a\in V(\mathcal{X}_1) $ ($ V(\mathcal{X}_2)$) considered while calculating the $SDI$ of $\mathcal{X}_1$($\mathcal{X}_2$) is also a SDS containing $a$ in $\mathcal{X}$. Also, by definition of join of FGs all the edges connecting the vertices of  $\mathcal{X}_1$ and    $\mathcal{X}_2$ are strong. Therefore, all sets of the form $\{u,v\}$ where $u\in (\mathcal{X}_1) $ and  $v\in (\mathcal{X}_2) $ are also SDSs in $\mathcal{X}$. Since the $_{sd}d(a)$ of a vertex $a$ is the minimum of weight of MSDSs containing $a$, $_{sd}d_{\mathcal{X}}(a)\leq_{sd}d_{\mathcal{X}_i}(a),i=1,2.$ Also, since the vertices of $\mathcal{X}$ is the union of vertices of $\mathcal{X}_1$ and $\mathcal{X}_2$, it is obtained that $SDI(\mathcal{X})\leq SDI(\mathcal{X}_1)+SDI(\mathcal{X}_2).$
\end{proof}
\section{Application}
\begin{figure}[h]
\centering
    \includegraphics[height=7cm,width=8cm]{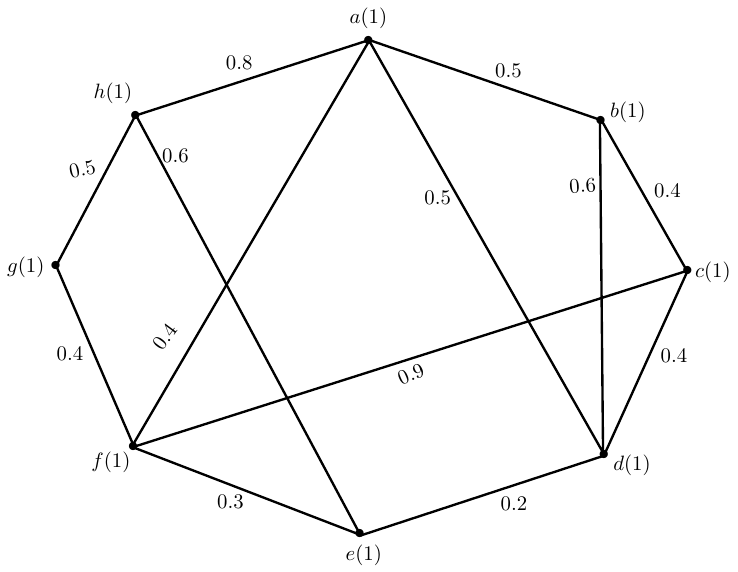}
    \caption{A FG model for road network within a city.}
    \label{fig10}
\end{figure}

This section discuss applications of SDD of a vertex. \\
Consider a city and the road network connecting different parts of the city. The model can be represented using an FG in which the vertices represent major parts of the city and edges represent the roads connecting them. Assign weight 1 to each vertex. The edge weight is the normalized number of accidents identified on each of the roads connecting these parts. The problem is to set up hospitals or clinics in the city so that almost all of the accidents can be addressed without any inconvenience. Suppose that while setting up these hospitals and clinics any one of the hospitals receives an innovative unit that consists of all the new treatment facilities. Then the idea of SDD can be used to decide in which part of the city the facilities must be provided so that more people get benefited from it. Consider FG in Figure \ref{fig10}.\\
Here a strong edge $ab$ means the number of accidents identified on direct road from $a$ to $b$ is greater than or equal to the number of accidents recorded on other alternate roads from $a$ to $b$ involving other parts of the city. Similarly, a $\delta-$ edge $ab $ means the number of accidents identified on direct road from $a$ to $b$ is less than the number of accidents recorded on other alternate roads. Hence $\delta -$ edges are neglected in this case. Now, consider a vertex $a,$ and let $\mathcal{D}$ be a minimal SDS containing $a$. Let the hospitals and clinics be set up at the parts of the city corresponding to the vertices in $\mathcal{D}.$ Also suppose that the extra facilities be provided at $a$. Then SDD of $a$ is the minimum normalised number of accidents that can be addressed by the hospitals and clinics covering the entire city with the extra facilities being provided at $a$. Now consider the SDDs of all vertices and take the vertex with maximum SDD say $v$ and let the minimal SDS considered for SDD of $v$ be $\mathcal{D}'$. Then it means maximum number of accidents are addressed if the hospitals or clinics are set up at parts corresponding to the vertices in $\mathcal{D}',$ and the extra facilities are provided at $v.$  For the FG in Figure \ref{fig10}, the $\delta-$ edges are $fe$ and $ed$. The vertex $e$ has the maximum SDD with the MSDS being $\{e,c,f\}$ and weight being 1.4. Hence the hospitals are to be set up at $e,c$ and $f$ with extra facilities at $e$. 
\\\\
Environmental pollution is a serious issue that must be addressed immediately. With the abrupt change in the public's perception of climate change, more organizations are becoming concerned about pollution and its impact on the environment. Long-term impacts of land filling include poor soil quality and contaminated water supplies, leading to crop failure, health risks, and other problems. Proper disposal of dry waste has several benefits. Recycling and composting reduce the amount of waste that ends up in landfills, and thus the negative impact on water and air quality. Recycling and composting are economically beneficial as they create jobs, boost the economy, and reduce the cost of waste disposal. Recycling also helps to conserve natural resources and minimizes the chemical residues that result from raw material extraction and processing.
The problem is to set up a waste drop-off or collection center in a municipality by minimizing the distance traveled by the collection agents or the people in the municipality. A FG is used to model the problem and the idea of SDD of a vertex is used to solve the problem. Consider a FG whose vertices represent wards in a municipality and whose edges represent the roads between the wards. Each vertex is assigned a weight of 1, and the edges receive a normalized weight based on the "nearness" concept. That means an edge of weight 0.9 between two wards $a$ and $b$ indicates that the wards $a$ and $b$ are very near. An edge of weight 0.1 between two wards $a'$ and $b'$ means wards $a'$ and $b'$ are very far. A strong edge between $a$ and $b$ indicates that the direct road from $a$ to $b$ is nearer than indirect road connecting $a$ to $b.$ Similarly a $\delta-$ edge between $a$ and $b$ indicates that the direct road from $a$ to $b$ is farther than indirect road connecting $a$ and $b.$ Hence in this problem the $\delta-$ edges are neglected. Now suppose for the waste management, the agents go door to door for waste collection and bring it to collection centres. People can also drop-off the dry wastes in centres allotted for the purpose. So, centers must be built so that the distance traveled by the agents and people is less. Also suppose that the population in a particular ward is very high. Then it is practical to build a centre in that particular ward. For this purpose the concept of MSDS is used. Consider a ward $a$ with the highest population. Then a MSDS $\mathcal{D}$, containing $a$ with the least weight serves a solution to the problem. If $\mathcal{D}=\{a,b,c\}$, is a  MSDS containing $a$ with least weight $w$, then that means the centres must be built in wards $a$ and $b$ and $c$ and the sum of maximum distance traveled to each centre is $w.$ A MSDS containing $a$ having maximum weight can also be considered to minimise the distance travelled by people or agents. Also consider all the SDDs and the vertex with maximum SDD can be considered for some special purpose as in the previous case so that the distance travelled is the least. 

\vspace{1cm}
\textbf{\large{Algorithm to find a minimal SDS containing a particular vertex}}

\vspace{1cm}

Consider a FG, $\mathcal{X}=(\mathcal{V},\varrho,\upsilon)$ with $n$ vertices, $v_1,v_2,...,v_n. $ Suppose it is required to find a minimal SDS containing a particular vertex say $v_k$. Then the following algorithm can be used to find a minimal SDS containing $v_k$. \\
First, a matrix is obtained from the FG which contains information about the strong edges and the weight of strong edges between the vertices. Consider an $n\times n $ matrix $A=[a_{ij}]$, where $$a_{ij} = \begin{cases}
 \upsilon(v_i,v_j)   &\text{ if $v_iv_j$ is a strong edge}\\
 0 & \text{ if $v_iv_j$ is a $\delta$ edge}
 \end{cases}$$
Now, consider $i^{th}$ row of the matrix $A$, then
$N_s(v_i)=\{v_{j}: \, a_{ij}>0, \, j=1,2,...n\}$. \\
$$N_s[v_i]=\{v_i\}\cup\{v_{j}: \, a_{ij}>0, \, j=1,2,...n\}$$
The problem is to find a minimal SDS containing a particular vertex say $v_{k_0}.$  
\\
Follow the following steps to obtain the required minimal SDS.
\begin{enumerate}
    \item Take vertex $v_{k_0}$. If $N_s[v_{k_0}]=\mathcal{V}$, then $\{v_{k_0}\}$ is a required set.
    \item Suppose $N_s[v_{k_0}]\neq \mathcal{V}$, take any vertex $v_{k_1}\neq v_{k_0}$ such that \\
    $$N_s[v_{k_0}]\setminus N_s[v_{k_1}] \neq \phi \text{ and } N_s[v_{k_1}]\setminus N_s[v_{k_0}] \neq \phi$$
    \item If $N_s[v_{k_1}]\cup N_s[v_{k_0}]=\mathcal{V}$, then $\{v_{k_1}, v_{k_0}\}$ is a required set. 
    \item Suppose $N_s[v_{k_1}]\cup N_s[v_{k_0}]\neq \mathcal{V}$, take next vertex $v_{k_2}\neq v_{k_0}, v_{k_1}$, such that \\
    $$N_s[v_{k_0}]\setminus\{ N_s[v_{k_1}] \cup N_s[v_{k_2}]\}\neq \phi \text{ and } N_s[v_{k_1}]\setminus \{N_s[v_{k_0}]\cup N_s[v_{k_1}]\} \neq \phi \text{ and } N_s[v_{k_2}]\setminus\{ N_s[v_{k_1}] \cup N_s[v_{k_0}]\}\neq \phi .$$
    \item If $N_s[v_{k_2}] \cup N_s[v_{k_1}]\cup N_s[v_{k_0}]=\mathcal{V}$, then $\{v_{k_2},v_{k_1},v_{k_0}\}$ is the required set.
    \item Continue the process until the whole vertex set is obtained. 
\end{enumerate}
Check all the possible combinations of vertices satisfying the conditions in each step of the algorithm. \\
After obtaining the minimal SDS, the weight of the set can also be found. \\
Suppose $D=\{v_{k_0}, v_{k_1}, v_{k_2},..., v_{k_m}\}$ is a minimal SDS obtained from the algorithm, then $$\mathcal{W}(D)= \sum \limits_{l=0}^m min\{a_{k_{\small{l}}j}: \, j=1,2,...,n\}.$$
After obtaining all the minimal SDSs containing $v_{k_0}$ check the weights of each set and take the set with least weight. The minimum of weight of minimal SDSs containing $v_{k_0}$ is the SDD of $v_{k_0}.$

\vspace{0.7cm}
\textbf{Illustration of Algorithm}
\vspace{0.7cm}

Consider the FG in Figure \ref{2}. In order to find a minimal SDS containing $a_4$ follow the following steps:
First obtain the matrix $A=[aij], \, 1\leq i,j,\leq n.$ From the Figure \ref{2}, matrix \\
$$A=\begin{bmatrix}
0 & 0 & 0 & 0.05 & 0 & 0 & 0 & 0\\
0 & 0 & 0 & 0.5 & 0 & 0 & 0 & 0\\
0 & 0 & 0 & 0.3 & 0 & 0 & 0 & 0\\
0.05 & 0.5 & 0.3 & 0 & 0.2 & 0.2 & 0 & 0\\
0 & 0 & 0 & 0.2 & 0 & 0.4 & 0.3 & 0\\
0 & 0 & 0 & 0.2 & 0.4 & 0 & 0 & 0.2\\
0 & 0 & 0 & 0 & 0.3 & 0 & 0 & 0\\
0 & 0 & 0 & 0 & 0 & 0.2 & 0 & 0\\
\end{bmatrix}$$
From the matrix, $N_s[a_1]=\{a_1,a_4\}$, $N_s[a_2]=\{a_2,a_4\}$, $N_s[a_3]=\{a_3,a_4\}$, $N_s[a_4]=\{a_4,a_1.a_2,a_3,a_5a_6\}$, $N_s[a_5]=\{a_5,a_4,a_6,a_7\}$, $N_s[a_6]=\{a_6,a_4,a_5,a_7\}$, $N_s[a_7]=\{a_7,a_5\}$ and $N_s[a_8]=\{a_8,a_6\}.$ 
\begin{enumerate}
    \item Here $N_s[a_4]\neq \mathcal{V}. $
    \item Take next vertex $a_1.$ But, $N_s[a_1]\setminus N_s[a_4]=\phi$.Therefore $a_1$ cannot be considered.
    \item By the same argument $a_2$ and $a_3$ is not considered.
    \item Now, take $a_5$. Here, $N_s[a_4]\setminus N_s[a_5]\neq \phi $ and $N_s[a_5]\setminus N_s[a_4]\neq \phi $. But $N_s[a_4]\cup N_s[a_5]\neq \mathcal{V}$.
    \item Hence consider next vertex $a_6$. Then, $N_s[a_4]\setminus \{N_s[a_5]\cup N_s[a_6]\}\neq \phi $, $N_s[a_5]\setminus \{N_s[a_4] \cup N_s[a_6]\}\neq \phi $ and $N_s[a_6]\setminus \{N_s[a_4] \cup N_s[a_5]\}\neq \phi $. Also, $N_s[a_4]\cup N_s[a_5]\cup N_s[a_6]=\mathcal{V}.$
    \item Hence $\{a_4,a_5,a_6\}$ is a required set with weight 0.45.
\end{enumerate}
Similarly, $\{a_4, a_7, a_8\}$ is also a set. Hence the given algorithm can be used find a minimal SDS containing a particular vertex of a FG. 
\section{Conclusion}
The topological index plays a vital role in various scientific fields, especially in the study of molecular structures and their properties. It provides valuable insights into the structural and chemical characteristics of molecules, enabling researchers to make predictions and understand fundamental properties without the need for extensive experimental data. The notion of domination in FGs has applications in various fields, including decision-making, pattern recognition, and image processing. It allows for the representation and analysis of uncertain or imprecise data, where the degrees of membership and connectivity can be fuzzy or uncertain. Here, both the ideas are combined to obtain the SDI of a FG. It is demonstrated that it is always possible to find a MSDS containing a particular vertex in a connected SFG. The concepts of upper strong domination number, strong irredundance number, strong upper irredundance number, strong independent domination number, and strong independence number are discussed and illustrated in the sections that follow. The SDD of each vertex is used to define the strong domination index (SDI) of a graph. The idea is used in a variety of fuzzy graphs, including complete, complete bipartite, complete r-partite FGs, fuzzy tree, fuzzy cycle, etc. The SDD and SDI-related inequalities are obtained. The research also takes into account the union and join of FG. Later sections offer applications for SDD of a vertex. An approach to obtain an MSDS containing a particular vertex is provided as an algorithm.  \\\\
\textbf{\large{Acknowledgement}}\\\\
The first author gratefully acknowledges the financial support of Council of Science and Industrial Research (CSIR), Government of India.\\
The authors would like to thank the DST, Government of India, for providing support to carry out this work under the scheme 'FIST' (No.SR/FST/MS-I/2019/40).

\end{document}